\documentclass[12pt,reqno]{amsart}
\usepackage[widespace]{fourier}
\usepackage{hyperref}
\hypersetup{colorlinks=true,linkcolor=blue, citecolor=blue}
\usepackage{amsmath,amscd,amsfonts,amssymb}
\usepackage[all]{xypic}
\usepackage[utf8]{inputenc}
\usepackage{longtable}
\numberwithin{equation}{section}
\newtheorem{theorem}{Theorem}[section]
\newtheorem{definition}[theorem]{Definition}

\newtheorem{lemma}[theorem]{Lemma}
\newtheorem{proposition}[theorem]{Proposition}

\newtheorem{remark}[theorem]{Remark}

\begin{document}

\title{Rationality of Moduli Space over reducible curve}

\author[A. Dey]{Arijit Dey}

\address{Department of Mathematics, Indian Institute of Technology-Madras, Chennai, India }

\email{arijitdey@gmail.com}

\author[Suhas, B N]{Suhas, B N}
\address{Department of Mathematics, Indian Institute of Technology-Madras, Chennai, India }
\email{chuchabn@gmail.com}

\maketitle
\abstract
Let $M(2,\textbf{\underline{w}},\chi)$ be the moduli space of rank $2$ torsion-free sheaves over a reducible nodal curve with each component having utmost two nodal singularities. We show that in each component of $M(2,\textbf{\underline{w}},\chi)$, the closure of rank $2$ vector bundles with fixed determinant is rational.
\endabstract
\section{Introduction}
Let $C$ be a connected projective curve over an algebraically closed field $k$ of characteristic $0$ having $N$ smooth components $C_i$ of genus $g_i \,\geq \, 2$ and $N-1$ nodes $P_i$ such that $C_i \cap C_{i+1}\,=\,P_i,$ for $i \,=\,1,2,\cdots,N-1$. We call such a curve as a {\it chain-like curve}.  Let $\textbf{\underline{w}}:\,=\,(w_1,\,w_2,\,\cdots,\,w_N)$ be an $N$-tuple of positive rational numbers such that $\sum_{j\,=\,1}^{N} \omega_j\,=\,1$, we call this a
{\it polarisation} on $C$. Let $\chi$ be an odd integer and $M(2,\textbf{\underline{w}},\chi)$ be the moduli space of rank $2$, ${\textbf{\underline{w}}}$-semi-stable torsion free sheaves with fixed Euler characteristic $\chi$ \cite{S2}. It is known that for a generic  choice of ${\textbf{\underline{w}}}$,  $M(2,\textbf{\underline{w}},\chi)$ has $2^{N-1}$
irreducible components $M_l(2,\textbf{\underline{w}},\chi)$, $l \,=\,1$ to $2^{N-1}$ \cite{T-B}. Each component is determined by choosing an N-tuple of integers $(\chi_1,\chi_2,\cdots,\chi_{N})$ satisfying inequalities \eqref{Inequality_for_chi_i} and \eqref{chi(E)} (for $n$\,=\,$2$)  such that for a generic vector bundle $E$ in the component, $\chi(E\mid_{C_i}) \,=\,\chi_i$.

Let $\xi$ be a line bundle on $C$ given by $(L_1,\,L_2,\, \cdots \,,\,L_N)$, where $L_i$'s are invertible sheaves on $C_i$'s respectively. Let
$\overline{M_l}(2,\textbf{\underline{w}},\chi,\xi)$ denote the closure of collection of vector bundles  with determinant
$\xi$ in $M_l(2,\textbf{\underline{w}},\chi)$. In this article we want to prove that this subvariety $\overline{M_l}(2,\textbf{\underline{w}},\chi,\xi)$ is rational for each $l$. When $N\,=\,2$ this result has appeared in \cite{BP-DA-S} and these subvarieties are interpreted as fixed determinant moduli space of torsion free sheaves \cite{BS}, \cite{N-S}. When $N \,>\, 2$ such an analogue of fixed determinant moduli space of torsion free sheaves is not known. We expect that if we have a similar notion of fixed determinant moduli space,  then our result will tell that it will be rational.

Over a smooth projective curve of genus $g \,\geq \, 2,$ the rationality of the moduli space was first proved by Tjurin \cite[Theorem 14]{TJ} in the rank $2$ and odd \text{\text{deg}}ree case. When rank and \text{\text{deg}}ree are coprime this result was generalized by Newstead \cite{N}, \cite{N2}, King and Schofield  \cite{K-S} in higher order of generalities.
It is still not known if the moduli space is rational or not in the non-coprime case, even for rank $2$ and \text{\text{deg}}ree $0$. In the non-smooth case, when the curve is
irreducible and has any number of nodal singularities and genus $\,\geq \, 2,$ rationality in the coprime case was proved by Bhosle and Biswas \cite[Theorem 3.7]{B-B}.  Over a
reducible nodal curve with two components  (i.e. $N\,=\,2$) the moduli space of torsion free sheaves with fixed determinant has two components \cite{BS}. The proof of
rationality of each of these components given in \cite{BP-DA-S} uses Nagaraj-Seshadri's description of the moduli space in terms of triples \cite{N-S}. For higher values of $N$,
such a description is not known. Hence the proof given in  \cite{BP-DA-S} can not be generalized for arbitrary $N$. The proof in this article is based on  Newstead's
idea \cite{N} and Teixidor I Bigas's description of the moduli space \cite{T-B} but involves several technical challenges. In fact Teixidor I Bigas's description of
the moduli space holds for more general curve known as {\it tree-like curve} but the combinatorics involved will be more complicated. It will be interesting as well as challenging to investigate rationality question in this case.

{\bf Acknowledgement.} We thank Prof. P. E. Newstead  for suggesting this problem and for a careful reading of the manuscript. We also thank Sarang Sane for some useful 
discussions.

\section{Description of the Moduli Space}
 Let $C$ be a chain-like curve with $N$ irreducible components $C_i$ of genus $g_i \,\geq \, 2$  such that $C_i \cap C_{i+1} \,=\, P_i,$ for $i \,=\,1,2,\cdots,N-1.$ Suppose $E$ is a vector bundle on $C$ of rank $n$ and Euler characteristic $\chi$ and $E_i$ is $E_{|_{C_i}}.$ Then one has the following exact sequence -
 \begin{equation}\label{Exact_Sequence_for_E}
 0\rightarrow E \rightarrow \bigoplus_{j\,=\,1}^{N}\alpha_{j_*}(E_j) \rightarrow T_{E} \rightarrow 0,
\end{equation}
where $\alpha_j:\, C_j \rightarrow C$ is a closed immersion for each $j$ and $T_E$ is a torsion sheaf supported only at the nodal points. Let $\chi_j:=\chi(E_j).$ Then from the exact sequence \eqref{Exact_Sequence_for_E}, one can show that
\begin{equation}\label{chi(E)}
 \chi\,= \,\sum_{j\,=\,1}^{N}\chi_j - n(N-1).
\end{equation}

Now let $\textbf{\underline{w}}:=(w_1,w_2,\cdots,w_N)$ be a polarization on $C,$ i.e., $w_j$ is a positive rational number for each $j,$ and $\sum_{j\,=\,1}^{N}w_j \,=\, 1.$ For any torsion-free sheaf $E$ on $C,$ we define
\begin{equation}\label{mu(E)}
 \mu(E)\,=\,\frac{\chi(E)}{\sum_{j\,=\,1}^{N}w_j~\text{rk}~(E_j)},
\end{equation}
where $E_j\,=\,\frac{E_{|_{C_j}}}{\text{torsion}}.$

\begin{definition}\label{Definition of stability}
 We say that a torsion-free sheaf  $E$ on $C$ is stable (resp. semi-stable) if for every proper sub-sheaf $G$ of $E,$ we have
 $\mu(G) \,<\, \mu(E)$ (resp. $\leq$).
\end{definition}

It is a theorem of Seshadri that over any reducible curve, the moduli space $M(n,\textbf{\underline{w}},\chi)$ of semi-stable torsion-free sheaves of rank $n$ and euler characteristic $\chi$ exists and is compact (see \cite[Chap VII]{S2}).

In the case of a chain-like curve $C$ with $N$ components, Teixidor i Bigas has proven in \cite[Theorem-1, Step-1]{T-B} that $M(n,\textbf{\underline{w}},\chi)$ has $n^{N-1}$ components for any generic polarization. \footnote{In fact she proves this result for tree-like curves.} In fact she has shown that if $E$ is a stable vector bundle on $C$ with Euler characteristic $\chi$ and $E_i$ has Euler characteristic $\chi_i$ for each $i,$ then $\chi_i$'s are going to satisfy the inequalities : \footnote{ The inequalities \eqref{Inequality_for_chi_i} follow from \cite[Theorem-1, Step-1, (1)]{T-B}.}
 \begin{equation}\label{Inequality_for_chi_i}
  (\sum_{j\,=\,1}^{i}w_j) \chi - \sum_{j\,=\,1}^{i-1}\chi_j + n(i-1) \,<\, \chi_i  \,<\,  (\sum_{j\,=\,1}^{i}w_j) \chi - \sum_{j\,=\,1}^{i-1}\chi_j + ni,
 \end{equation}
for $i\,=\,1,2,\cdots,N-1,$ provided $(\sum_{j\,=\,1}^{i}w_j) \chi$ is not an integer for each $i \in \{1,2,\cdots,N-1\}.$ She also proves in \cite[Theorem-1, Step-2]{T-B} that, for any choice of a semi-stable vector bundle $E_i$ on each component $C_i$ with Euler characteristic $\chi_i$ as in the inequality \eqref{Inequality_for_chi_i}, and gluing by any isomorphism at the nodes, one can obtain a semi-stable vector bundle $E$ on $C$ and further if one of the $E_i$ is stable, so is $E.$ Since there are $n^{N-1}$ choices for such tuples $(\chi_1,\chi_2,\cdots,\chi_N),$ one can conclude that $M(n,\textbf{\underline{w}},\chi)$ has $n^{N-1}$ components, each component corresponding to a particular type of $(\chi_1,\chi_2,\cdots,\chi_N).$

In what follows, we assume that $\chi$ is odd and $n\,=\,2.$ We also choose the polarization $\textbf{\underline{w}}$ in such a way that stability coincides with semi-stability.

\section{Construction of a stable family}
%Throughout this section, we fix $\chi=1.$ %Suppose $L_j$ is an invertible sheaf on $C_j$ for each $j.$
Let $\chi\,=\,1$ and $\chi_1,\chi_2,\cdots,\chi_{N-1}$ be integers satisfying the inequalities
\begin{equation}\label{Inequality_for_chi_i_corresponding_to_chi_equal_to_one}
 (\sum_{j\,=\,1}^{i}w_j) - \sum_{j\,=\,1}^{i-1}\chi_j + 2(i-1) \,<\,  \chi_i \,<\,  (\sum_{j\,=\,1}^{i}w_j) - \sum_{j\,=\,1}^{i-1}\chi_j + 2i.
\end{equation}
Let $\chi_N$ be an integer which fits into the following equation-
\begin{equation}\label{chi}
 \chi\,=\,1\,=\,\sum_{j\,=\,1}^{N}\chi_j - 2(N-1).
\end{equation}
So there are $2^{N-1}$ choices for the tuple $(\chi_1,\chi_2,\cdots,\chi_N)$ satisfying \ref{Inequality_for_chi_i_corresponding_to_chi_equal_to_one} and \ref{chi}. By the equation \eqref{chi}, for all these choices one has
\begin{equation}\label{sum-of_chi_i}
\sum_{j\,=\,1}^{N}\chi_j \,=\, 2N-1.
\end{equation}
Let $L_j$ be an invertible sheaf on $C_j$ for each $j \in \{1,2,\cdots,N\}.$
\begin{definition}
 We say that the tuple $(L_1,L_2,\cdots,L_N)$ is of type $(\chi_1,\chi_2,\cdots,\chi_N)$ if $deg~(L_j)\,=\,\chi_j-2(1-g_j)$ for each $j.$
\end{definition}
Throughout this section we fix an invertible sheaf $L_j$ on $C_j$ for each $j$, such that $L_j$'s are globally generated and the tuple $(L_1,L_2,\cdots,L_N)$ is of type $(\chi_1,\chi_2,\cdots,\chi_N)$ where $\chi_1,\chi_2,\cdots,\chi_{N-1}$ are as in \ref{Inequality_for_chi_i_corresponding_to_chi_equal_to_one} and $\chi_N$ is as in the equation \eqref{chi}.

%Now let $P_1$ be the node on $C_1,$ $P_{N-1}$ be the node on $C_{N}$ and $P_{j-1},P_{j}$ be the nodes on $C_{j},$ for $j=2,\cdots,N-1.$
%Now suppose $P_{j_1},P_{j_2}, \cdots,P_{j_l}$ are the nodes on $C_{j},$ where $l$ depends on $j.$ For each $k \in \{1,2,\cdots,l\},$
Let
\begin{equation*}
 T_{j} \,=\, \{t \in H^{0}(C_j,L_j)|~ t(P_{j-1}) \,\neq \, 0~\text{and}~t(P_j) \,\neq \, 0 \},
\end{equation*}
for $j\,=\,2,\cdots,N-1.$ Similarly let
\begin{equation*}
 T_{1} \,=\, \{t \in H^{0}(C_1,L_1)|~ t(P_{1}) \,\neq \, 0 \},
\end{equation*}
and
\begin{equation*}
 T_{N} \,=\, \{t \in H^{0}(C_N,L_N)|~ t(P_{N-1}) \,\neq \, 0 \}.
\end{equation*}
Clearly for each $j \in \{1,\cdots,N\},$ $T_{j}$ is a non-empty Zariski-open subset of the affine space $H^{0}(C_j,L_j).$ So there are sections in $H^{0}(C_j,L_j)$ which do not vanish on any node of $C_j.$ Let $s_j \in H^{0}(C_j,L_j)$ be one such section for each $j.$

%Let $Z=\{ (i,j)|~i,j\in \{1,2,\cdots,N\},~i <j ~and~ C_i \cap C_j \neq \emptyset\}.$ Clearly cardinality of $Z$ is $N-1,$ which is the number of nodes on $C.$ For each $(i,j) \in Z,$ let $C_i \cap C_j= \{P_{ij}\}.$ Now let
Let
\begin{equation*}
 \lambda_{j} :\, L_j(P_{j}) \rightarrow L_{j+1}(P_{j})
\end{equation*}
be the linear map of one dimensional vector spaces which sends $s_j(P_{j})$ to $s_{j+1}(P_{j}),$ where $j\,=\,1,2,\cdots,N-1.$ We now define an invertible sheaf $L$ on $X$ as follows:
 \begin{equation*}
  L \,=\, \{(t_1,t_2,\cdots,t_N) \in \bigoplus_{j\,=\,1}^{N}\alpha_{j_*}(L_j)~|~ t_j(P_{j})\xrightarrow{\lambda_{j}} t_{j+1}(P_{j})\}.
 \end{equation*}

Clearly $(s_1,s_2,\cdots,s_N) \in H^{0}(C,L)$ where $s_j$'s are as defined above. From now on, we call this section as the "distinguished section". By definition of $L$ we have the following short exact sequence:
\begin{equation} \label{Equation 3.1}
 0\rightarrow L \rightarrow \bigoplus_{j\,=\,1}^{N}\alpha_{j_*}(L_j) \rightarrow T \rightarrow 0,
\end{equation}
where $T$ is the torsion sheaf which is supported at the nodal points. Since there are $(N-1)$ nodes, $H^{0}(C,T)$ is a vector space of \text{\text{\text{\text{dim}}}}ension $N-1.$

\begin{remark}
 Let $L \,=\,(L_1,\cdots,L_N)$ be as above. The fact that $\chi_j$'s are chosen as in \ref{Inequality_for_chi_i_corresponding_to_chi_equal_to_one} and \eqref{chi} will imply that $deg~(L_j) \, \geq \,  2g_j-1$ for each $j.$ If for some $j$ $\text{deg}~(L_j) \,=\, 2g_j-1,$ then such an $L_j$ need not be globally generated in general. But one can always choose an invertible sheaf $L_j$ of degree $2g_j-1$ such that $L_j$ is globally generated (see \cite[Remark 3.2(a)]{BP-DA-S}).
\end{remark}

\begin{lemma}\label{Lemma 3.2}
 Let $L$ be as above. Then

 (i) The functor $H^{0}(C,-)$ applied to \eqref{Equation 3.1} is exact.

  (ii) \text{\text{\text{\text{dim}}}}~($H^{0}(C,L))=g$ and \text{\text{\text{\text{dim}}}} ($H^{1}(C,L)) \,=\, 0.$

  (iii) \text{\text{\text{\text{dim}}}}~($H^{0}(C,L^{*})) \,=\, 0$ and \text{\text{\text{\text{dim}}}} ($H^{1}(C,L^{*})) \,=\, 3g-2,$ where $L^{*}$ is the dual of $L.$

 \begin{proof}
  Applying the functor $H^{0}(C,-)$ to \eqref{Equation 3.1}, we get the exact sequence

  \begin{equation}\label{Equation 3.2}
   0\rightarrow H^{0}(C,L) \rightarrow \bigoplus_{j\,=\,1}^{N}H^{0}(C_j,L_j) \xrightarrow{\beta} H^{0}(C,T).
  \end{equation}
  We claim $\beta$ is surjective. Consider the set  $ D \,=\, \{(s_1,0,\cdots,0),(0,s_2,0,\cdots,0),\cdots,(0,0,\cdots,0,s_{N-1},0)\},$ where $s_j$'s are the components of the "distinguished section". Clearly this is a linearly independent set in $\bigoplus_{j\,=\,1}^{N}H^{0}(C_j,L_j).$ Since the "distinguished section" goes to zero under $\beta,$ it is clear that image of each element of $D$ under $\beta$ is non-zero and in fact $\beta(D)$ is linearly independent in $H^{0}(C,T).$ Hence $\beta$ is surjective. This proves $(i).$

  Now by the choice of $L_{j}$'s it is clear that \text{\text{\text{\text{dim}}}}~ ($\bigoplus_{j\,=\,1}^{N}H^{0}(C_j,L_j)) \,=\, g + (N-1)$ and \text{\text{\text{\text{dim}}}}~ ($\bigoplus_{j \,=\,1}^{N}H^{1}(C_j,L_j)) \,=\,0$. So by $(i)$ and \eqref{Equation 3.2}, \text{\text{\text{\text{dim}}}}~ ($H^{0}(C,L)) \,=\, g$ and by taking the long exact sequence associated to \eqref{Equation 3.1}, we can conclude that \text{\text{\text{\text{dim}}}}~ ($H^{1}(C,L)) \,=\,0.$ This proves $(ii).$

  To prove $(iii),$ again by the choice of $L_j$'s, it is clear that $\text{deg}~(L_j^*) \,<\, 0$ for each $j.$ So $H^0(C_j,L_j^*) \,=\, 0$ for each $j.$ Since $L^* \hookrightarrow \bigoplus_{j\,=\,1}^{N}\alpha_{j_*}(L_j^*),$ it is clear that $H^0(C,L^*) \,=\,0.$ So

  \begin{eqnarray*}
   \text{\text{\text{\text{dim}}}}~ (H^{1}(C,L^{*})) &=& - \chi(L^{*}) \\
                        &=& -(\chi(\bigoplus_{j=1}^{N}\alpha_{j_*}(L_j^*))-(N-1)) \\
                        &=& -((\sum_{j=1}^{N}(\text{deg}~(L_j^*)+(1-g_j)))-(N-1)) \\
                        &=& 3g-2.
  \end{eqnarray*}
 \end{proof}
\end{lemma}

Now let $\gamma$ be a proper subset of $\{1,2,\cdots,N\}.$ Let
\begin{equation}\label{Definition of V-gamma}
 V_{\gamma} \,=\, \{(t_1,t_2,\cdots,t_N) \in H^{0}(C,L)~|~ t_{i} \neq 0,~ \text{if i}~ \in \gamma ~ \text{and}~ t_{i}=0 ~ \text{otherwise}\}.
\end{equation}
Clearly $\overline{V_{\gamma}}$, the closure of $V_{\gamma}$ in $H^{0}(C,L)$, is a linear subspace of $H^{0}(C,L)$ and hence is closed and irreducible in Zariski topology.
\begin{lemma}\label{Lemma 1.3}
 Let $\gamma \,=\, \{i_1,i_2,\cdots,i_t\}$ be a proper subset of $\{1,2,\cdots,N\}$ such that $i_1,i_2,\cdots,i_t$ are consecutive integers. If $i_1 \,=\, 1$ or $i_t \,=\, N,$
 \begin{equation*}
  \text{\text{\text{\text{dim}}}}~(\overline{V_{\gamma}}) \, \leq \, \sum_{j\,=\,1}^{t}h^{0}(C_{i_j},L_{i_j})-t.
 \end{equation*}
 Otherwise
 \begin{equation*}
  \text{\text{\text{\text{dim}}}}~(\overline{V_{\gamma}}) \, \leq \, \sum_{j\,=\,1}^{t}h^{0}(C_{i_j},L_{i_j})-(t+1).
 \end{equation*}

 \begin{proof}
  If $i_1 \,=\, 1$ or $i_t\,=\, N,$ the union $C_{i_{1}} \cup \cdots \cup C_{i_{t}}$ has $t-1$ internal nodes and one external node. If $i_1 \,\neq \, 1$ and $i_t \,\neq \, N,$ the union $C_{i_{1}} \cup \cdots \cup C_{i_{t}}$ has $t-1$ internal nodes and two external nodes.  So by the definition of $\overline{V_{\gamma}},$ the Lemma follows.
 \end{proof}
\end{lemma}
Let
\begin{equation}\label{Definition of V}
V = \bigcup_{\gamma}\overline{V_{\gamma}}.
\end{equation}
Since the "distinguished section" belongs to $ H^{0}(C,L) \backslash V,$ we can conclude that $V$ is a proper closed subset of the affine space $H^{0}(C,L).$ Let
\begin{equation}\label{Compliment of V}
R = H^{0}(C,L) \backslash V.
\end{equation}
 Clearly every element of $R$ defines an injective map $\mathcal{O}_C \hookrightarrow L$ and conversely if any non-zero section $(t_1,t_2,\cdots,t_N)$ of $L$ defines an injective map $\mathcal{O}_C \hookrightarrow L,$ then such a section should be in $R,$ for otherwise, it belongs to $V$ which means $t_i \,=\, 0$ for some $i$, and so such a section $(t_1,t_2,\cdots,t_N)$ cannot define an injective map $\mathcal{O}_C \hookrightarrow L.$ So we have
\begin{equation}\label{R}
 R \,=\, \{\psi \in H^{0}(C,L)~|~ \psi:\mathcal{O}_{C} \hookrightarrow L ~ \text{is injective}\}.
\end{equation}

\begin{lemma}\label{Dimension V-gamma}
 Let $\gamma \,=\, \{i_1,i_2,\cdots,i_t\}$ be a proper subset of $\{1,2,\cdots,N\}$ such that $i_1 \,<\, i_2 \,<\, \cdots \,<\, i_t.$ Then
\begin{equation*}
 \text{\text{\text{\text{dim}}}}~(\overline{V_{\gamma}}) \,\leq \, \sum_{j\,=\,1}^{t}g_{i_{j}}.
\end{equation*}
 As a consequence,
\begin{equation*}
 \text{\text{\text{\text{dim}}}}~(V_{\gamma}) \,\leq \, \sum_{j\,=\,1}^{t}g_{i_{j}}.
\end{equation*}

\begin{proof}
 By Reimann-Roch theorem and the choice of the invertible sheaves $L_j,$ we know that
 \begin{eqnarray}\label{Dimension-H0}
  \sum_{j\,=\,1}^{t}h^{0}(C_{i_{j}},L_{i_{j}}) &=& \sum_{j\,=\,1}^{t}[\text{\text{\text{\text{deg}}}}~(L_{i_{j}}) + (1-g_{i_{j}})] \nonumber \\
                                           &=& \sum_{j\,=\,1}^{t}[\chi_{i_{j}}-2(1-g_{i_{j}}) + (1-g_{i_{j}})] \nonumber \\
                                           &=& \sum_{j\,=\,1}^{t}\chi_{i_{j}} - t + \sum_{j=1}^{t} g_{i_{j}}.
 \end{eqnarray}
 With this in mind, we prove the Lemma by considering two different cases for $\gamma.$

 \textbf{Case~A:}
  We first assume that $\gamma$ consists of consecutive integers. We now consider three sub-cases -

 \textbf{Case~i:} Suppose $i_1\,=\,1.$ This implies $i_t \,=\, t.$ Now since the $\chi_j$'s satisfy \eqref{Inequality_for_chi_i_corresponding_to_chi_equal_to_one},
  $\sum_{j\,=\,1}^{t}\chi_j$ is either equal to $2t$ or equal to $2t-1.$ In any case,
 \begin{equation}\label{chi-Case(i)}
  \sum_{j\,=\,1}^{t}\chi_j \,\leq \, 2t.
 \end{equation}
 Therefore by the equation \eqref{Dimension-H0} and the inequality \eqref{chi-Case(i)},
 \begin{eqnarray}\label{h0-Case(i)}
  \sum_{j\,=\,1}^{t}h^{0}(C_{j},L_{j}) & \leq & 2t -t + \sum_{j\,=\,1}^{t} g_{j} \nonumber \\
                                           &=& \sum_{j\,=\,1}^{t} g_{j} + t.
 \end{eqnarray}
 Combining the inequality \ref{h0-Case(i)} and the Lemma \ref{Lemma 1.3}, we have
 \begin{equation*}
  \text{\text{\text{\text{dim}}}}~(\overline{V_{\gamma}}) \,\leq \, \sum_{j\,=\,1}^{t}g_j,
 \end{equation*}
 which proves the Lemma for this case.

 \textbf{Case~ii:} Suppose $i_t \,=\,N.$ This implies $i_1-1 \,=\, N-t.$ Now from the equation
 \begin{equation*}
  \sum_{i\,=\,1}^{N}\chi_i \,=\, \sum_{i\,=\,1}^{i_{1}-1}\chi_i + \sum_{i\,=\,i_{1}}^{N}\chi_i,
 \end{equation*}
 we have
 \begin{equation}\label{chi-Case(ii)}
  \sum_{i\,=\,i_{1}}^{N}\chi_i \,=\, \sum_{i\,=\,1}^{N}\chi_i - \sum_{i\,=\,1}^{i_{1}-1}\chi_i.
 \end{equation}
 The choice of $\chi_i$'s will imply that the sum $\sum_{i\,=\,1}^{i_{1}-1}\chi_i$ has to be at least $2i_1-3.$ Using this fact and the equation \eqref{sum-of_chi_i} in the equation \eqref{chi-Case(ii)}, we get
\begin{eqnarray}\label{Chi-Case(ii)}
 \sum_{i\,=\,i_{1}}^{N}\chi_i & \leq & (2N-1) - (2i_{1}-3) \nonumber \\
                          & = & (2N-1) - 2(N-t+1)+3 \nonumber \\
                          & = & 2t.
\end{eqnarray}
Combining this with the equation \eqref{Dimension-H0}, we get
 \begin{eqnarray}\label{h0-Case(ii)}
  \sum_{j\,=\,1}^{t}h^{0}(C_{i_{j}},L_{i_{j}})  & \leq & 2t -t + \sum_{j\,=\,1}^{t} g_{i_{j}} \nonumber \\
                                           &=& \sum_{j\,=\,1}^{t} g_{i_{j}} + t.
 \end{eqnarray}
  From the inequality \ref{h0-Case(ii)} and the Lemma \ref{Lemma 1.3}, we have
 \begin{equation*}
  \text{\text{\text{\text{dim}}}}~(\overline{V_{\gamma}}) \,\leq \, \sum_{j\,=\,1}^{t}g_{i_{j}}.
 \end{equation*}

 \textbf{Case~iii:} Suppose $i_1 \,\neq \, 1$ and $i_t \,\neq \, N.$ Then since
 \begin{equation*}
  \sum_{i\,=\,1}^{i_t}\chi_i \,=\, \sum_{i\,=\,1}^{i_1-1}\chi_i + \sum_{i\,=\,i_1}^{i_t}\chi_i,
 \end{equation*}
 we have
 \begin{eqnarray}\label{chi-Case(iii)}
  \sum_{i\,=\,i_1}^{i_t}\chi_i & = & \sum_{i\,=\,1}^{i_t}\chi_i - \sum_{i\,=\,1}^{i_1-1}\chi_i \nonumber \\
                           & \leq & 2i_t - (2i_1-3) \nonumber \\
                           & = & 2(i_1+(t-1)) - (2i_1-3) \nonumber \\
                           & = & 2t+1.
 \end{eqnarray}
 So combining with the equation \eqref{Dimension-H0}, we get
 \begin{equation}\label{h0-Case(iii)}
  \sum_{j\,=\,1}^{t}h^{0}(C_{i_{j}},L_{i_{j}}) \,\leq \, \sum_{j\,=\,1}^{t} g_{i_{j}} + (t+1).
 \end{equation}
 So by the inequality \ref{h0-Case(iii)} and the Lemma \ref{Lemma 1.3}, we can conclude that
\begin{equation*}
  \text{\text{\text{\text{dim}}}}~(\overline{V_{\gamma}}) \,\leq \, \sum_{j\,=\,1}^{t}g_{i_{j}}.
 \end{equation*}
 This proves the Lemma for \textbf{Case~A}.

 \textbf{Case~B:}
 Now suppose $\gamma$ is such that $i_1,i_2,\cdots,i_t$ are not consecutive.
 Let $\mathcal{C}_1,\cdots,\mathcal{C}_r$ be the connected components of $C_{i_{1}} \cup \cdots \cup C_{i_{t}}.$ Clearly each $\mathcal{C}_i$ consists of either single irreducible component or some consecutive irreducible components of $C.$ Let the corresponding subset of indices be $\gamma_i.$ Then $\gamma_i$ is either singleton or has consecutive integers. Therefore
 \begin{eqnarray}\label{Dimension V-gamma-Case(B)}
  \text{\text{\text{\text{dim}}}}~(\overline{V_{\gamma}}) & = & \sum_{\gamma_i \equiv \mathcal{C}_i}\text{\text{\text{\text{dim}}}}~(\overline{V_{\gamma_i}}) \nonumber \\
   & \leq & \sum_{i_j \in \gamma}g_{i_{j}}.
 \end{eqnarray}
 The last inequality comes by \textbf{Case~A}.
 This proves the Lemma.
\end{proof}
\end{lemma}

\begin{remark}\label{Remark 3.3}
  Let $q \,=\, 3g-2.$ Then by fixing a basis of $H^{1}(C,L^{*}),$ we can identify it with $k^{q}.$ We have the natural $k^{*}-$action on $k^{q}$ and
      \begin{equation*}
      W \,=\, \{ (a_{1},a_{2}, \cdots,a_{q}) \in k^{q}~| ~a_{1} \,\neq \, 0 \}
      \end{equation*}
   is clearly an invariant Zariski-open subset of $k^{q}$ under this $k^{*}-$ action.

   Let $A:= \{(a_{1},a_{2}, \cdots,a_{q})\in W ~|~ a_{1} \,=\,1 \}$ (Clearly $A$ is Zariski closed  and every orbit of $k^{*}-action$ on $W$ meets $A$ in exactly one point).
\end{remark}

\begin{lemma}(cf.~\cite{H-L}, \cite{N-R}, \cite{N})\label{Lemma 3.4}
 Let $L$ be as above. Then there exists a vector space $V^{\prime}$ and a universal extension
 \begin{equation}\label{Equation 3.3}
  0 \rightarrow \mathcal{O}_{C \times V^{\prime}} \rightarrow \tilde{\mathcal{E}} \rightarrow \pi^{*}(L) \rightarrow 0
 \end{equation}
of bundles over $C \times V^{\prime}$ (where $\pi:C \times V^{\prime} \rightarrow C$ is the projection map), such that there is a natural isomorphism
\begin{equation*}
 \alpha:V^{\prime} \rightarrow H^{1}(C,L^{*})
\end{equation*}
where for each $v \in V^{\prime}$, $\alpha(v)$ is the element corresponding to the restriction of the extension \eqref{Equation 3.3} to $\{v\} \times C.$
\end{lemma}

\begin{remark}\label{Remark 3.5}
 Suppose $\tilde{\mathcal{E}}$ is as in Lemma \ref{Lemma 3.4} and $v \in H^{1}(C,L^{*})$ is such that $\text{\text{\text{\text{dim}}}}~(H^{0}(C,\tilde{\mathcal{E}}_{v}))\,=\,1.$ Then one can easily see that for any $w \in H^{1}(C,L^{*}),$ $\tilde{\mathcal{E}}_{v} \,\cong \, \tilde{\mathcal{E}}_{w}$ if and only if $v$ and $w$ are in the same orbit under the natural action of $k^{*}$ on $H^{1}(C,L^{*}).$
\end{remark}

\begin{lemma}\label{Lemma 3.6}
  Let $L$ be as above. Then there exists an extension
\begin{equation}\label{Equation 3.4}
 0\rightarrow \mathcal{O}_{C} \rightarrow E \rightarrow L \rightarrow 0 ,
\end{equation}
for which $\text{\text{\text{\text{dim}}}}~ (H^{0}(C,E))\,=\,1$ and this extension can be chosen to correspond to a point of $A,$ where $A$ is as in Remark \ref{Remark 3.3} above.

\begin{proof}
 Suppose $a \in H^{1}(C,L^{*})$ and \eqref{Equation 3.4} is the corresponding extension. For any section $0 \neq \delta \in H^{0}(C,L),$ we have a non-trivial morphism
  \begin{eqnarray}\label{Equation 3.5}
    \delta:\,\mathcal{O}_{C} & \longrightarrow & L .
  \end{eqnarray}
  Tensoring \eqref{Equation 3.5} by the dualizing sheaf $\omega_{C}$ and applying the global section functor, we get the map
 \begin{equation*}
  H^{0}(C,\omega_{C}) \longrightarrow H^{0}(C,L \otimes \omega_{C}).
 \end{equation*}
Taking dual and using the duality theorem, we get the map
 \begin{equation*}
  H^{1}(C,L^{*}) \xrightarrow{\tilde{\delta}}  H^{1}(C,\mathcal{O}_{C}) .
 \end{equation*}
 This implies
 \begin{eqnarray}
  \text{\text{\text{\text{dim}}}}~(\text{ker}~(\tilde{\delta})) & \geq & \text{\text{\text{\text{dim}}}}~(H^{1}(C,L^{*})) -g \nonumber ~\,>\, ~0.
 \end{eqnarray}
 Applying the sheaf functors $\mathcal{H}om(L,-)$ and $\mathcal{H}om(\mathcal{O}_{C},-)$ to \eqref{Equation 3.4} and taking the long exact sequence, we get the following
 commutative diagram -
$$\xymatrix{
0 \ar@{->}[r]^{} \ar@{->}[d]^{} & Hom(L,\mathcal{O}_{C}) \ar@{->}[r]^{} \ar@{->}[d]^{} & Hom(L,E) \ar@{->}[r]^{} \ar@{->}[d]^{} & Hom(L,L) \ar@{->}[r]^{} \ar@{->}[d]^{} & H^{1}(C,L^{*}) \ar@{->}[r]^{}  \ar@{->}[d]^{\tilde{\delta}}  & \cdots \\
0 \ar@{->}[r]^{} & H^{0}(C,\mathcal{O}_{C})  \ar@{->}[r]^{} & H^{0}(C,E)  \ar@{->}[r]^{} & H^{0}(C,L) \ar@{->}[r]^{} & H^{1}(C,\mathcal{O}_{C}) \ar@{->}[r]^{} & \cdots
}$$
From this diagram, it is clear that $\delta$ lifts to a section on $E$ if and only if $\tilde{\delta}(a)\,=\,0.$ This fact is proved in \cite[Lemma 3.1]{N-R}, in greater generality. Also
\begin{eqnarray*}
 \text{\text{\text{\text{dim}}}}~(H^{0}(C,E)) &=& \text{\text{\text{\text{dim}}}}~(H^{0}(C,\mathcal{O}_{C})) + \text{\text{\text{\text{dim}}}}~(\text{ker}~(H^{0}(C,L)) \rightarrow H^{1}(C,\mathcal{O}_{C})).
\end{eqnarray*}
(One can prove that $\text{\text{\text{\text{dim}}}}~(H^{0}(C,\mathcal{O}_{C}))\,=\,1$ and $\text{\text{\text{\text{dim}}}}~(H^{1}(C,\mathcal{O}_{C}))\,=\,g$ by using the arguments similar to those in Lemma \ref{Lemma 3.2}).
So
\begin{eqnarray*}
 \text{\text{\text{\text{dim}}}}~(H^{0}(C,E))~ \,>\,~1 & \Leftrightarrow & \text{\text{\text{\text{dim}}}}~(\text{ker}~(H^{0}(C,L)) \rightarrow H^{1}(C,\mathcal{O}_{C})) \,\geq \, 1 \\
                       & \Leftrightarrow & a \in p_{1}(Y),
\end{eqnarray*}
where $Y$ is a subset of $H^{1}(C,L^{*}) \times H^{0}(C,L)$ defined by $(a,\delta) \in Y \Leftrightarrow 0 \,\neq \, \delta$ and $\tilde{\delta}(a)\,=\,0,$ and $p_{1}$ is the first projection.

We want an extension of the form \eqref{Equation 3.4} such that $\text{\text{\text{\text{dim}}}}~(H^{0}(C,E))\,=\,1.$ It is enough to show that $\text{\text{\text{\text{dim}}}}~(p_{1}(Y)) \,\leq \, \text{\text{\text{\text{dim}}}}~(H^{1}(C,L^{*})) -1.$

Let $Y^{\prime} \subseteq Y$ be such that $(a,\delta) \in Y^{\prime} \Leftrightarrow \delta: \, \mathcal{O}_{C} \hookrightarrow L~ \text{and} ~ \tilde{\delta}(a)\,=\,0.$ Since $C$ is a stable curve  \cite[Definition I.I]{D-M}, the dualizing sheaf $\omega_{C}$ is locally free \cite[Theorem 1.2]{D-M}. So if
\begin{equation*}
  \delta: \, \mathcal{O}_{C} \hookrightarrow L ,
\end{equation*}
is an injective map, then tensoring by $\omega_{C},$ we get an injective map
\begin{equation*}
 \omega_{C} \hookrightarrow L \otimes \omega_{C}.
\end{equation*}
So the induced map
\begin{equation*}
 H^{1}(C,L^{*}) \xrightarrow{\tilde{\delta}}  H^{1}(C,\mathcal{O}_{C})
\end{equation*}
will be surjective and
\begin{equation}\label{Equation 3.6}
  \text{\text{\text{\text{dim}}}}~(\text{ker}~(\tilde{\delta}))\,=\, \text{\text{\text{\text{dim}}}}~(H^{1}(C,L^{*}))-g ~\,>\,~ 0.
\end{equation}
It is clear from \eqref{Equation 3.6} that $Y^{\prime}$ is non-empty. We claim that $Y^{\prime}$ is open in $Y.$

%Let $R^{\prime} = H^{0}(C,L) \backslash R,$ where $R$ is the open set defined above (see equations \eqref{Compliment of V}, \eqref{R}).
Consider the second projection
 \begin{equation*}
   p_{2}:\, H^{1}(C,L^{*}) \times H^{0}(C,L) \rightarrow H^{0}(C,L).
 \end{equation*}
 As $p_{2}^{-1}(V)$ is closed in $H^{1}(C,L^{*}) \times H^{0}(C,L),$~ $p_{2}^{-1}(V) \cap Y$ is closed in $Y,$ where $V$ is as defined in \eqref{Definition of V}.

  Clearly $p_{2}^{-1}(V) \cap Y \,=\, Y \backslash Y^{\prime}.$
  Therefore $Y^{\prime}$ is open in $Y.$  We claim  $\text{\text{\text{\text{dim}}}}~(Y^{\prime})\,=\,\text{\text{\text{\text{dim}}}}~(Y).$

  If $Y$ is irreducible, or if every irreducible component of $Y$ intersects $Y^{\prime},$ then we are done. Otherwise let $Y_1 \subset Y$ be an irreducible component of $Y$ such that $Y_1 \cap Y^{\prime} \,=\,\emptyset.$ Consider the map $p_{2_{|_{Y_1}}}: \, Y_1 \rightarrow \overline{p_2(Y_1)},$ where $\overline{p_2(Y_1)}$ is the closure of $p_2(Y_1)$ in $H^{0}(C,L).$ Since $Y_1$ is irreducible, $\overline{p_2(Y_1)}$ is also irreducible, and being a closed sub-variety of the affine space $H^{0}(C,L),$ it is an affine variety. Let $\mathcal{U} \subset Y_1$ be an affine open and consider the map $p_{2_{|_{\mathcal{U}}}}:\mathcal{U} \rightarrow \overline{p_2(Y_1)}.$ This map is clearly a dominant map of irreducible affine varieties. Hence there exists an open subset $W_1$ of $\overline{p_2(Y_1)}$ such that $W_1 \subset p_2(Y_1).$ Now further restricting $p_2$ to $p_2^{-1}(W_1)$ we get a surjective map $p_2^{-1}(W_1) \rightarrow W_1$ of irreducible varieties. Therefore by \cite[Theorem 1.25 (ii),Page-75]{Sh}, there exists an open subset $W_2$ in $W_1$ such that
  \begin{equation}\label{fiber dimension}
   \text{\text{\text{\text{dim}}}}~(p_2^{-1}(\phi)) \,=\, \text{\text{\text{\text{dim}}}}~(p_2^{-1}(W_1)) - \text{\text{\text{\text{dim}}}}~(W_1),
  \end{equation}
  for all $\phi \in W_2.$ But it is clear that $\text{\text{\text{\text{dim}}}}~(p_2^{-1}(W_1))\,=\, \text{\text{\text{\text{dim}}}}~(Y_1)$ and $\text{\text{\text{\text{dim}}}}~(W_1) \,=\, \text{\text{\text{\text{dim}}}}~(\overline{p_2(Y_1)}).$ So the equation \eqref{fiber dimension} becomes
  \begin{equation*}
   \text{\text{\text{\text{dim}}}}~(p_2^{-1}(\phi))\,=\, \text{\text{\text{\text{dim}}}}~(Y_1) - \text{\text{\text{\text{dim}}}}~(\overline{p_2(Y_1)}),
  \end{equation*}
  for all $\phi \in W_2.$
 So we have
 \begin{equation}\label{dimension Y1}
  \text{\text{\text{\text{dim}}}}~(Y_1) \,=\, \text{\text{\text{\text{dim}}}}~(p_2^{-1}(\phi)) + \text{\text{\text{\text{dim}}}}~(\overline{p_2(Y_1)}),
 \end{equation}
for all $\phi \in W_2.$ (In equations \eqref{fiber dimension} and \eqref{dimension Y1}, by $~p_2^{-1}(\phi)$ we mean $~p_2^{-1}(\phi) \cap Y $).
We know that $p_2^{-1}(\phi) \cap Y \,=\, \text{ker}~\tilde{\phi}.$ So by equation \eqref{dimension Y1}, to find $\text{\text{\text{\text{dim}}}}~(Y_1)$ (or at least an "optimal" upper bound for $\text{\text{\text{\text{dim}}}}~(Y_1)$), we have to find $\text{\text{\text{\text{dim}}}}~(\overline{p_2(Y_1)})$ and $\text{\text{\text{\text{dim}}}}~(\text{ker}~\tilde{\phi})$ for some $\phi \in W_2.$

  Now since $Y_1 \cap Y^{\prime} \,=\,\emptyset,$ $p_2(Y_1) \subset V,$ where $V$ is as defined in the equation \eqref{Definition of V}. So there exists a $\gamma$ such that $\overline{p_2(Y_1)} \subset \overline{V_{\gamma}}.$ Let $\gamma{\prime}$ be a proper subset of $\{1,2,\cdots,N\}$ such that
  \begin{equation*}
   \text{\text{\text{\text{dim}}}}~(\overline{V_{\gamma^{\prime}}})\,=\, \text{min}~\{\text{\text{\text{\text{dim}}}}~(\overline{V_{\gamma}})~|~p_2(Y_1) \subset \overline{V_{\gamma}}\}.
  \end{equation*}
  Then $p_2(Y_1)$ has to intersect $V_{\gamma^{\prime}}$ itself, for otherwise, $p_2(Y_1)$ will be completely inside a smaller dimensional $\overline{V_{\gamma}}$ which contradicts the fact that $\overline{V_{\gamma^{\prime}}}$ is minimum dimensional among all such $\overline{V_{\gamma}}.$ Since $V_{\gamma^{\prime}}$ is open in its closure, $V_{\gamma^{\prime}} \cap \overline{p_2(Y_1)}$ is open in $\overline{p_2(Y_1)}.$ So it has to intersect $W_2$ and $W_2 \cap V_{\gamma^{\prime}} \cap \overline{p_2(Y_1)}$ is non-empty, open and is a subset of $p_2(Y_1).$ We denote this open set by $W_3.$ Clearly $\text{\text{\text{\text{dim}}}}~(\overline{p_2(Y_1)})\,=\, \text{\text{\text{\text{dim}}}}~(W_3).$ So the equation \eqref{dimension Y1} becomes -
  \begin{equation}\label{Dimension-Y1}
   \text{\text{\text{\text{dim}}}}~(Y_1) \,=\, \text{\text{\text{\text{dim}}}}~(p_2^{-1}(\phi)) + \text{\text{\text{\text{dim}}}}~(W_3),
  \end{equation}
  where $\phi \in W_3.$

  Let $\phi \,=\, (t_1,t_2,\cdots,t_N) \in W_3$ be arbitrary. Since $W_3 \subset V_{\gamma^{\prime}},$ $t_i\,=\,0$ for $i$ not in $\gamma^{\prime}$ and $t_i \,\neq \, 0$ for $i \in \gamma^{\prime}.$ So for $i \in \gamma^{\prime},$ $t_i: \, \mathcal{O}_{C_i} \rightarrow L_i$ is injective and the induced map $\tilde{t_i}:\, H^{1}(C,L^{*}) \rightarrow H^{1}(C,\mathcal{O}_{C})$ is surjective. Now consider the following commutative diagram -
  $$\xymatrix{
 H^{1}(C,L^{*}) \ar@{->}[r]^{\tilde{\phi}} \ar@{->}[d] & H^{1}(C,\mathcal{O}_C) \ar@{->}[d]
\\ \bigoplus_{j\,=\,1}^{N} H^{1}(C_j,L_j^{*}) \ar@{->}[r]^{(\tilde{t_1},\cdots,\tilde{t_N})} & \bigoplus_{j\,=\,1}^{N} H^{1}(C_j,\mathcal{O}_{C_j}).
}$$
Both the vertical arrows in the above diagram are surjective because they are gotten by taking the long exact sequence corresponding to the short exact sequences -
\begin{equation*}
 0\rightarrow L^{*} \rightarrow \bigoplus_{j\,=\,1}^{N}\alpha_{j_{*}}(L_{j}^{*}) \rightarrow T^{\prime} \rightarrow 0,
\end{equation*}
and
\begin{equation*}
 0\rightarrow \mathcal{O}_C \rightarrow \bigoplus_{j\,=\,1}^{N}\alpha_{j_{*}}(\mathcal{O}_{C_j}) \rightarrow \tilde{T} \rightarrow 0
\end{equation*}
respectively and observing that $T^{\prime}$ and $\tilde{T}$ are supported only at the nodal points $P_1,\cdots,P_{N-1},$ and hence $H^{1}(C,T^{\prime})\,=\, 0\,=\,H^{1}(C,\tilde{T}).$ Since $\tilde{t_i}$ is surjective for $i \in \gamma{\prime},$ we have $\text{\text{\text{\text{dim}}}}~(\text{Im}~(\tilde{t_i})) \,=\, \text{\text{\text{\text{dim}}}}~(H^{1}(C_i,\mathcal{O}_{C_i})\,=\,g_i,$ for such $i,$ and $\text{\text{\text{\text{dim}}}}~(\text{Im}~(\tilde{t_i})) \,=\, 0,$ for $i$ not in $\gamma^{\prime}.$ Now since the above diagram commutes, we have $\text{\text{\text{\text{dim}}}}~(\text{Im}~(\tilde{\phi})) \,\geq \, \sum_{i \in \gamma^{\prime}}g_i.$ This implies
\begin{equation}\label{Kernel phi-tilde}
 \text{\text{\text{\text{dim}}}}~(\text{Ker}~(\tilde{\phi})) \,\leq \, \text{\text{\text{\text{dim}}}}~(H^{1}(C,L^{*}))- \sum_{i \in \gamma^{\prime}}g_i ,
\end{equation}
for all $\phi \in W_3.$ Also since $W_3 \subset V_{\gamma^{\prime}},$ we have $\text{\text{\text{\text{dim}}}}~(W_3) \,\leq \, \sum_{i \in \gamma^{\prime}}g_i$ by Lemma \ref{Dimension V-gamma}. Using this and the inequality \eqref{Kernel phi-tilde} in the equation \eqref{Dimension-Y1}, we get
\begin{equation*}
 \text{\text{\text{\text{dim}}}}~(Y_1) \,\leq \, \text{\text{\text{\text{dim}}}}~(H^{1}(C,L^{*})).
\end{equation*}
So if $Y_1$ is an irreducible component of $Y$ such that $Y^{\prime} \cap Y_1 \,=\, \emptyset,$ then by the above arguments
\begin{equation}\label{Final Expression for Dimension-Y1}
 \text{\text{\text{\text{dim}}}}~(Y_1) \,\leq \, \text{\text{\text{\text{dim}}}}~(H^{1}(C,L^{*})).
\end{equation}

  Now by using similar arguments as above and the equation \eqref{Equation 3.6}, one can prove
  \begin{eqnarray*}
   \text{\text{\text{\text{dim}}}}~(Y^{\prime}) &=& \text{\text{\text{\text{dim}}}}~(H^{0}(C,L))+ \text{\text{\text{\text{dim}}}}~(H^{1}(C,L^{*})) -g \\
                  &=& \text{\text{\text{\text{dim}}}}~(H^{1}(C,L^{*})).
  \end{eqnarray*}
  This implies by \eqref{Final Expression for Dimension-Y1} that $\text{\text{\text{\text{dim}}}}~(\text{Y}) \,=\, \text{\text{\text{\text{dim}}}}~(H^{1}(C,L^{*})).$ This proves our claim that $\text{\text{\text{\text{dim}}}}~(\text{Y})\,=\, \text{\text{\text{\text{dim}}}}~(Y^{\prime}) \,=\, \text{\text{\text{\text{dim}}}}~(H^{1}(C,L^{*})).$

  Now, if $a \in p_{1}(Y),$ then $\text{\text{\text{\text{dim}}}}~(p_{1}^{-1}(a) \cap Y) \,\geq \, 1$ since $(a,\delta) \in Y \Leftrightarrow (a,\lambda \delta) \in Y$ for every non-zero scalar $\lambda.$ So
  \begin{eqnarray*}
   \text{\text{\text{\text{dim}}}}~(p_{1}(Y)) & \leq & \text{\text{\text{\text{dim}}}}~(Y) -1 \\
                &=& \text{\text{\text{\text{dim}}}}~(H^{1}(C,L^{*})) -1 .
  \end{eqnarray*}
  This proves the first part of the lemma.

  Now let
  \begin{equation*}
   B \,=\, \{ a \in H^{1}(C,L^{*}) ~|~ \text{\text{\text{\text{dim}}}}~(H^{0}(C,\mathcal{\tilde{E}}_{a}))\,=\,1\},
  \end{equation*}
   where $\mathcal{\tilde{E}}$ is as in Lemma \ref{Lemma 3.4}. Then by semi-continuity theorem, $B$ is clearly a $k^{*}-$ invariant open subset of $H^{1}(C,L^{*}).$ Let $W$ be as in Remark \ref{Remark 3.3} above. Then $B \cap W \neq \emptyset$ as both $B$ and $W$ are Zariski-open subsets of an affine space. Since $B$ and $W$ are $k^{*}-$ invariant, $B \cap W$ is also $k^{*}-$ invariant. Therefore $B$ meets $A$ also.
\end{proof}
\end{lemma}

 \begin{remark}\label{Remark 3.8}
(a) Let $L$ be as above and $F \subset L$ be a sub-sheaf. Now let $t \,\geq \, 1$ be an integer and $\gamma \,=\, \{i_1,i_2,\cdots,i_t\} \subset \{1,2,\cdots,N\}$ be a proper subset such that $i_1 \,<\, i_2 \,<\, \cdots, \,<\, i_t.$ Suppose $F$ is such that $rk(F_{i_1}),rk(F_{i_2}),\cdots,rk(F_{i_t})$ are all equal to one and $rk(F_i) \,=\,0$ for $i \,\neq \, i_1,i_2,\cdots,i_t.$ Then it is clear that $H^0(C,F) \subset \overline{V_{\gamma}},$ where $V_{\gamma}$ is as defined in equation \eqref{Definition of V-gamma}. So we have
 \begin{eqnarray}\label{h0(C,F)}
  h^0(C,F) & \leq & \text{dim}~(\overline{V_{\gamma}}) \nonumber \\
           & \leq & \sum_{j\,=\,1}^t g_{i_j}.
 \end{eqnarray}
 The last inequality follows from Lemma \ref{Dimension V-gamma}. Since $\chi(F) \leq h^0(C,F),$ from the inequality \ref{h0(C,F)} we have
 \begin{equation}\label{chi(F)}
  \chi(F) \,\leq \, \sum_{j\,=\,1}^t g_{i_j}.
 \end{equation}

(b) Let $E$ be a locally free sheaf on $C$ of rank two such that $\text{\text{\text{\text{dim}}}}~(H^{0}(C,E))\,=\,1.$ Let $G \subset E$ be a subsheaf such that its Euler characteristic is positive. Then $\text{\text{\text{\text{dim}}}}~(H^{0}(C,G)) \,=\,1.$ Moreover, if $E$ is an extension as in \eqref{Equation 3.4}, then the map $\mathcal{O}_{C} \rightarrow E$ factors through $G$ because the sections of $E$ are same as sections of $G.$
\end{remark}

\begin{lemma}\label{Lemma 3.9}
 Let $L$ be as above and $E$ be an extension as in \eqref{Equation 3.4} such that $\text{\text{\text{\text{dim}}}}~(H^{0}(C,E)) \,=\,1.$ Then $E$ is  stable.
 \begin{proof}
  Let $G \subset E$ be a proper sub-sheaf. Since the weights are chosen in such a way that semi-stability coincides with stability, it is enough to prove
  \begin{equation}\label{Equation 3.14}
   \chi(G) \,\leq \, (\sum_{j\,=\,1}^{N}w_j ~ \text{rk}~(G_j))\frac{\chi(E)}{2}.
  \end{equation}
   By the choice of $E$ and $L$ it is clear that $\chi(E)\,=\,1.$ So we have to prove
  \begin{equation*}
   \chi(G) \,\leq \, \frac{(\sum_{j\,=\,1}^{N}w_j~ \text{rk}~(G_j))}{2}.
  \end{equation*}
   We prove this by considering all possible cases for $G.$

  $\mathbf{Case~ 1:}$ Suppose $G$ is such that $\text{rk}~(G_i)\,=\,2$ for all $i.$ So we have to prove $\chi(G) \,\leq \, 1$ in this case.

   Suppose $\chi(G) \,>\, 1.$ Then $\text{\text{\text{\text{dim}}}}~(H^{0}(C,G)) \,>\, 1 \,=\, \text{\text{\text{\text{dim}}}}~(H^{0}(C,E)).$ But this is not possible as $G \subset E.$ So we are done.

 $\mathbf{Case~ 2:}$ Suppose $\text{rk}~(G_j) \,=\,1$ for all $j.$ In this case we have to prove $\chi(G) \,\leq \, \frac{1}{2}.$

 Suppose $\chi(G) \,>\, \frac{1}{2},$ then by Remark \ref{Remark 3.8}, $\text{\text{\text{\text{dim}}}}~(H^{0}(C,G))\,=\,1$ and the map $\mathcal{O}_{C} \hookrightarrow E$ factors through $G.$ Let us denote $(\frac{G}{\mathcal{O}_{C}})$ by $F$ for notational convenience. Since $F$ is a subsheaf of $(\frac{E}{\mathcal{O}_{C}})\,=\,L,$ it implies $F$ is either torsion free or zero. But the fact that $\text{rk}~(G_j)\,=\,1$ for all $j$ forces $F$ to be zero, for otherwise, $F$ is supported at finitely many points and so it cannot be a subsheaf of $L.$ This implies $G \cong \mathcal{O}_C$ which gives a contradiction to our assumption that $\chi(G) \,>\, \frac{1}{2}.$

 $\mathbf{Case~ 3:}$ Suppose $\text{rk}~(G_j)\,=\,0$ for some $j.$ We want to prove $\chi(G) \,\leq \, \frac{(\sum_{j\,=\,1}^{N}w_j~ \text{rk}~(G_j))}{2}.$  Suppose $\chi(G) \,>\, \frac{(\sum_{j\,=\,1}^{N}w_j~ \text{rk}~(G_j))}{2}.$ This implies $\chi(G)$ is positive and so by the arguments in the Remark \ref{Remark 3.8},  $\text{\text{\text{\text{dim}}}}~(H^{0}(C,G))\,=\,1$ and the map $\mathcal{O}_{C} \hookrightarrow E$ factors through $G.$ So we have $\mathcal{O}_{C,p} \hookrightarrow G_p$ for each $p \in C.$ In particular $\mathcal{O}_{C,p} \hookrightarrow G_p$ for each smooth point $p \in C_j.$ This contradicts the fact that $\text{rk}~(G_j)\,=\,0$ for some $j.$

 $\mathbf{Case~ 4:}$ Suppose $i_1, i_2,\cdots,i_t$ are indices in the increasing order such that $\text{rk}~(G_{i_1}),\text{rk}~(G_{i_2}),\cdots,\text{rk}~(G_{i_t})$ are all equal to two and $\text{rk}~(G_{j})\,=\,1$ for $j \,\neq \, i_1,i_2,\cdots,i_t$ (here we are assuming that $G_j$'s are of mixed rank and $\text{rk}~(G_j) \,\neq \, 0$ for all $j$). Again, to prove the required result, if we assume on the contrary that $\chi(G) \,>\, \frac{(\sum_{j\,=\,1}^{N}w_j~ \text{rk}~(G_j))}{2},$ then the same arguments as before say that the map $\mathcal{O}_{C} \hookrightarrow E$ factors through $G.$ Denoting $(\frac{G}{\mathcal{O}_C})$ by $F$ as before, and using the inequality \ref{chi(F)} in the equation $\chi(G) \,=\, \chi(\mathcal{O}_C) + \chi(F),$ we get
 \begin{eqnarray}\label{Equation 3.17}
  \chi(G) & \leq & \chi(\mathcal{O}_C) + \sum_{j\,=\,1}^{t}g_{i_j} \nonumber \\
          & = & 1-g + \sum_{j\,=\,1}^{t}g_{i_j} \nonumber \\
          & = & 1- \sum_{k \,\neq \, i_j}g_k.
 \end{eqnarray}
 Since each $g_k \,\geq \, 2,$ the inequality \ref{Equation 3.17} implies that $\chi(G)$ is negative. This is a contradiction. So we are done.
 \end{proof}
\end{lemma}

\begin{remark}\label{Remark 3.10}
(a) In this section we have assumed $\chi\,=\,1$ and $\chi_i$'s satisfy \eqref{Inequality_for_chi_i_corresponding_to_chi_equal_to_one} and \ref{chi}. So, for example, all the results of this section are valid if $L$ is of type $(1,2,2,\cdots,2)$ or $(1,3,1,3,\cdots,1,3,1)$ or $(1,3,1,2,2,\cdots,2)$ and so on.

(b) Suppose $L$ is of type $(\chi_1,\chi_2,\cdots,\chi_N)$ and
   \begin{equation*}
    0\rightarrow \mathcal{O}_{C} \rightarrow E \rightarrow L \rightarrow 0
   \end{equation*}
 is an exact sequence, then $\chi_j$'s are precisely the Euler characteristics of $E_j$'s. Since $\chi(E) \,=\, \sum_{j\,=\,1}^{N}\chi_j - 2(N-1),$ if $\chi(E)$ is odd, then $\sum_{j\,=\,1}^{N}\chi_j$ should be odd. This implies that the cardinality of the set
 \begin{equation}\label{Equation 3.18}
  X \,=\, \{\chi_j~|~j \in \{1,2,\cdots,N\}~ \text{and} ~\chi_j~\text{is an odd integer}\}
 \end{equation}
 is an odd number.
\end{remark}

In the proof of Lemma \ref{Lemma 3.6}, we saw that the set $B \,=\, \{a \in H^{1}(C,L^{*}) ~|~ \text{\text{\text{\text{dim}}}}~(H^{0}(C,\tilde{\mathcal{E}}_{a}))\,=\,1 \}$ is a $k^{*}-$ invariant nonempty open subset of $H^{1}(C,L^{*})$ and $B \cap A$ is a non-empty open subset of the affine space $A.$ Let $S \,=\, B \cap A.$ Now, if $a \in B,$ then by Lemma \ref{Lemma 3.9}, it is clear that $\tilde{\mathcal{E}}_{a}$ is stable. So $S$ is a non-empty open subset  of the affine space $A$ consisting of stable rank two  locally free sheaves $\tilde{\mathcal{E}}_{a}$ such that $\text{\text{\text{\text{dim}}}}~(H^{0}(C,\tilde{\mathcal{E}}_{a}))\,=\,1.$ Since $S \subset A$ is Zariski-open, $\text{\text{\text{\text{dim}}}}~(S)\,=\,\text{\text{\text{\text{dim}}}}~(A)\,=\,3g-3.$

\begin{lemma}\label{Lemma 3.11}
 Let $L$ be an invertible sheaf on $C$ of any type mentioned above. Then there exists a non-empty open subset $S$ of an affine space and a locally free sheaf $\mathcal{E}^{\prime}$ of rank two on $S \times C$ such that

 (i) $\text{\text{\text{\text{dim}}}}~(S)\,=\,3g-3,$

 (ii) $\mathcal{E}^{\prime}_{s} \cong \mathcal{E}^{\prime}_{t}$  $\Leftrightarrow$  $s\,=\,t,$

 (iii) for all $s \in S,$ $\mathcal{E}^{\prime}_{s}$ is stable and $\Lambda^{2}(\mathcal{E}^{\prime}_{s})\,=\,L.$

 \begin{proof}
  Consider the sheaf $\tilde{\mathcal{E}}$ on $H^1(C,L^*) \times C$ obtained in \eqref{Equation 3.3} and restrict it to $S \times C,$ where $S$ is as defined just above. Let $\mathcal{E}^{\prime}\,=\,\tilde{\mathcal{E}}_{|_{S \times C}}.$ We have already seen that $\text{\text{\text{\text{dim}}}}~(S)\,=\,3g-3.$ By the definition of $S$ (more precisely, by the definition of $A$) and Remark \ref{Remark 3.5}, it is clear that $\mathcal{E}^{\prime}_{s} \cong \mathcal{E}^{\prime}_{t}$  $\Leftrightarrow$  $s\,=\,t.$ Again by the definition of $S$ it is clear that for all $s \in S,$ $\mathcal{E}^{\prime}_{s}$ is stable and $\Lambda^{2}(\mathcal{E}^{\prime}_{s})\,=\,L.$ This proves the Lemma.
 \end{proof}
\end{lemma}

\section{Rationality}
 Suppose $\chi\,=\,1,$ $C,$ $\mathbf{\underline{w}},$ $L$ and $\chi_j$'s are as before. Let $M_{(\chi_1,\cdots,\chi_N)}(L)$ denote the collection of all vector bundles in $M(2,\mathbf{\underline{w}},\chi\,=\,1)$ with determinant $L.$ Since $M(2,\mathbf{\underline{w}},\chi\,=\,1)$ is a coarse moduli space, by Lemma \ref{Lemma 3.11}, we have an injective morphism
\begin{equation}\label{Equation 4.2}
  f:\, S \rightarrow M(2,\mathbf{\underline{w}},\chi\,=\,1),
 \end{equation}
where $S$ is as in the Lemma \ref{Lemma 3.11}. But the image of $f$ lands in $M_{(\chi_1,\cdots,\chi_N)}(L).$ Since $S$ and $M_{(\chi_1,\cdots,\chi_N)}(L)$ are both smooth varieties of same \text{\text{\text{\text{dim}}}}ension and we are in characteristic zero, it implies that $f$ is birational. So $M_{(\chi_1,\cdots,\chi_N)}(L)$ is rational.

\begin{remark}\label{Remark 4.1}
Let $U^{\prime} \subset M_{(\chi_1,\cdots,\chi_N)}(L)$ be the open subset consisting of all vector bundles $E$ such that $E_i$ is semi-stable for each $i.$ By \cite[Step 2]{T-B}, $U^{\prime} \,\neq \, \emptyset.$ Let $U \,=\, f^{-1}(U^{\prime}).$ Then $U \,\neq \, \emptyset$ as $f$ is birational. Since $U \subset S$ is open, \text{\text{\text{\text{dim}}}}~(U)\,=\,3g-3.

Let $\mathcal{E}^{\prime}$ be as in Lemma \ref{Lemma 3.11}. Then $\mathcal{E}^{\prime}_{|_U}$ is a locally free sheaf of rank two on $U \times C$ such that for each $u \in U,$ $\mathcal{E}^{\prime}_u$ is stable and for each $i,$ $\mathcal{E}^{\prime}_{u}|_{C_i}$ is semi-stable. Clearly the restriction map
\begin{equation}\label{Equation 4.3}
  f:\, U \rightarrow M_{(\chi_1,\cdots,\chi_N)}(L)
 \end{equation}
 is birational.

 This remark is important because, stability of an arbitrary vector bundle $E$ on $C$ does not guarantee the semi-stability of $E_i$ on $C_i.$ But there is a non-empty open set in the moduli space consisting of stable vector bundles whose restriction to each component is semi-stable.
\end{remark}

With these in mind, we now state and prove the main proposition-

\begin{proposition}\label{Proposition 4.1}
  Let $\chi$ be an odd integer and $C, \mathbf{\underline{w}}$ be as mentioned above. Let $\chi_1, \cdots, \chi_{N-1}$ be the integers satisfying the  inequalities -
  \begin{equation*}
 (\sum_{j\,=\,1}^{i}w_j) \chi - \sum_{j\,=\,1}^{i-1}\chi_j + 2(i-1) \,<\, \chi_i  \,<\,  (\sum_{j\,=\,1}^{i}w_j) \chi - \sum_{j\,=\,1}^{i-1}\chi_j + 2i,
\end{equation*}
and $\chi_N$ be the integer such that
\begin{equation*}
 \chi \,=\, \sum_{j\,=\,1}^{N}\chi_j - 2(N-1).
\end{equation*}
Let $L \,=\,(L_1,\cdots,L_N)$ be an invertible sheaf on $C$ of type $(\chi_1,\chi_2,\cdots,\chi_N).$ Then there exists a non-empty open subset $U$ of an affine space and a locally free sheaf $\mathcal{E}$ on $U \times C$ of rank two such that

 (i) $\text{\text{\text{\text{dim}}}}~(U)\,=\,3g-3,$

 (ii) $\mathcal{E}_{u} \cong \mathcal{E}_{t}$  $\Leftrightarrow$  $u\,=\,t,$

 (iii) for all $u \in U,$ $\mathcal{E}_{u}$ is stable and $\Lambda^{2}(\mathcal{E}_{u})\,=\,L.$

 \begin{proof}
  Let $i_1,i_2,\cdots,i_t$ be the indices in the increasing order such that $\chi_{i_1},\chi_{i_2},\cdots,\chi_{i_t}$ are odd integers and $\chi_j$'s are even integers if $j \,\neq \, i_1,\cdots,i_t$. So $t$ is an odd number (see Remark \ref{Remark 3.10}(b) and equation \eqref{Equation 3.18}).

  If $\chi\,=\,1$ and each $L_j$ is globally generated, then by Lemma \ref{Lemma 3.11} and the above arguments, we are done. Also when $\chi\,=\,1,$ among all the $2^{N-1}$ choices for the tuple $(\chi_1,\cdots,\chi_N),$ there exists a choice for which $\chi_{i_j}\,=\,1$ if $j$ is odd, $\chi_{i_j}\,=\,3$ if $j$ is even and $\chi_i\,=\,2$ if $i$ does not belong to $\{i_1,i_2,\cdots,i_t\}.$ Now since $L\,=\,(L_1,\cdots,L_N)$ is of type $(\chi_1,\cdots,\chi_N),$ and $\chi_{i_j}$'s are odd for $j\,=\, 1, \cdots,t,$ $\text{\text{\text{\text{deg}}}}~(L_{i_j})\,=\, \chi_{i_j}-2(1-g_{i_j})$ will be odd for each $i_j$ and $\text{\text{\text{\text{deg}}}}~(L_{r})\,=\, \chi_r - 2(1-g_r)$ will be even for $r \,\neq \, i_j.$
  For notational convenience, we write $\text{\text{\text{\text{deg}}}}~(L_{i_j})\,=\,2l_{i_j}-1$ for $j\,=\,1,\cdots,t,$ and
  $\text{\text{\text{\text{deg}}}}~(L_r)\,=\,2l_r$ for $r \,\neq \, i_j,$ where $l_{i_j}\,=\, \frac{\chi_{i_j}-2(1-g_{i_j})+1}{2}$ for each $i_j,$ and $l_r\,=\,\frac{\chi_r - 2(1-g_r)}{2}$ for each $r \,\neq \, i_j.$  Let $M$ be an invertible sheaf on $C$ such that

  $\text{\text{\text{\text{deg}}}}~(M_r)\,=\, (l_r-g_r)$ for $r \,\neq \, i_2, i_4,\cdots,i_{t-1},$

  $\text{\text{\text{\text{deg}}}}~(M_{r})\,=\,(l_{r}-g_{r}-1)$ for $r \in \{i_2, i_4,\cdots,i_{t-1}\},$ and

  $L_r \otimes M_r^{-2}$ is globally generated for $r\,=\,1,2,\cdots,N$ (see \cite[Remark 4.2]{BP-DA-S} for the existence of such an $M$).

 It is clear that

 $\text{\text{\text{\text{deg}}}}~(L_r \otimes M_r^{-2})\,=\, 2g_r$ for $r \,\neq \, i_1,i_2,\cdots,i_t,$

 $\text{\text{\text{\text{deg}}}}~(L_{r} \otimes M_{r}^{-2})\,=\, 2g_{r}-1$ for $r \in \{i_1,i_3,\cdots,i_t\},$ and

 $\text{\text{\text{\text{deg}}}}~(L_{r} \otimes M_{r}^{-2})\,=\, 2g_{r}+1$ for $r \in \{i_2,i_4,\cdots,i_{t-1}\}.$

 So by Lemma \ref{Lemma 3.11} and Remark \ref{Remark 4.1}, there exists a non-empty open subset $U$ of an affine space and  a locally free sheaf $\mathcal{E}^{\prime}$ on $U \times C$ such that

 (i) $\text{\text{\text{\text{dim}}}}~(U) \,=\,3g-3,$

 (ii) $\mathcal{E}^{\prime}_{u} \cong \mathcal{E}^{\prime}_{t}$  $\Leftrightarrow$  $u\,=\,t,$

 (iii) for all $u \in U,$ $\mathcal{E}^{\prime}_{u}$ is stable and $\Lambda^{2}(\mathcal{E}^{\prime}_{u})\,=\,L \otimes M^{-2}.$

 Let $\mathcal{E}\,=\,\mathcal{E}^{\prime} \otimes p_C^*(M).$ Then $\mathcal{E}_u\,=\,\mathcal{E}^{\prime}_u \otimes M$ and $\Lambda^2(\mathcal{E}_u)\,=\,\Lambda^2(\mathcal{E}^{\prime}_u) \otimes M^{2}\,=\,L$ for all $u \in U.$

 Now since $\mathcal{E}^{\prime}_u $ is stable, by the Remark \ref{Remark 4.1}, $\mathcal{E}^{\prime}_{u}|_{C_i}$ is semi-stable for all $i.$ This means $\mathcal{E}^{\prime}_{u}|_{C_i} \otimes M_i$ is semi-stable for all $i$ and $\mathcal{E}^{\prime}_{u}|_{C_j} \otimes M_j$ is stable for $j \,=\, i_1,i_2,\cdots,i_t.$
 This proves that $\mathcal{E}_u$ is stable (see \cite[Step 2]{T-B}).

 This proves the proposition.
 \end{proof}
\end{proposition}
From the Proposition, we can conclude that the sheaf $\mathcal{E}$ on $U \times C$ induces a morphism
 \begin{equation*}
  f:\, U \rightarrow M_{(\chi_1,\cdots,\chi_N)}({L}).
 \end{equation*}
 By $(ii)$ of the Proposition, $f$ is injective. Since $U$ and $M_{(\chi_1,\cdots,\chi_N)}({L})$ are smooth varieties of same \text{\text{\text{\text{dim}}}}ension and we are in characteristic zero, it implies that $f$ is birational. So $M_{(\chi_1,\cdots,\chi_N)}({L})$ is rational.


\begin{thebibliography}{99}

%\bibitem{BDS} Barik, P;  Dey, A; Suhas, B. N. ;  On the rationality of Nagaraj-Seshadri moduli space, to appear in Bull. Sci. Math.

\bibitem{BP-DA-S} Barik, Pabitra; Dey, Arijit; Suhas, B.N.; \emph{On the Rationality of Nagaraj-Seshadri Moduli Space},~~~ Accepted for publication in Bulletin Des Sciences Mathematiques.

\bibitem{BS} Basu, Suratno; \emph{On a relative Mumford-Newstead theorem};~~~` to appear in Bull. Sci. Math.

\bibitem{B-B}  Bhosle, Usha N. ; Biswas, Indranil; Brauer group and birational type of moduli spaces of torsionfree sheaves on a nodal curve. Comm. Algebra 42 (2014), no. 4, 1769--1784.

\bibitem{K-S}  King, Alastair ; Schofield, Aidan; \emph{Rationality of moduli of vector bundles on curves},~~~Indag. Math. (N.S.) 10 (1999), no. 4, 519--535.

\bibitem{H-L} Lange , Herbert; \emph{Universal Families of Extensions}, Journal of Algebra 83, 101-112( 1983).

\bibitem{N-S} Nagaraj, D.S. ; Seshadri, C.S. ; \emph{\text{\text{Deg}}enerations of moduli spaces of vector bundles on curves -I}, ~~~ Proc. Indian Acad. Sci (Math. Sci.),Vol-107, No:2, 1997.

\bibitem{N-R} Narasimhan, M.S. ; Ramanan, S; \emph{Moduli of vector bundles on a compact Riemann surface},~~~Ann. of Math.89,14-51(1969).

\bibitem{N} Newstead, P.E. ; \emph{Rationality of Moduli Spaces of Stable Bundles},~~~Math.Ann. 215, 251-268(1975).

\bibitem{N2}  Newstead, P. E.; \emph{Correction to: "Rationality of moduli spaces of stable bundles"},~~~Math. Ann. 249, (1980), no. 3 , 281-282.

\bibitem{D-M} Deligne, P; Mumford, D; \emph{The irreducibility of space of curves of given genus},  Publications Mathematiques de l'I.H.E.S.,tome 36(1969),p.75-109.

\bibitem{T-B} Teixidor i Bigas, M ; \emph{Moduli spaces of (semi)stable vector bundles on tree-like curves}, Math.Ann. 290, 341-348 (1991).

\bibitem{Sh} Shafarevich , Igor R.; \emph{Basic Algebraic Geometry 1},~~~ Third edition, Springer.

%\bibitem{BP-DA-S} Barik, Pabitra; Dey, Arijit; Suhas, B.N.; \emph{On the Rationality of Nagaraj-Seshadri Moduli Space},~~~ Accepted for publication in Bulletin Des Sciences Mathematiques.

\bibitem{S2} Seshadri, C. S.; \emph{Fibre vectoriels sur les courbes algebrriques},~~~Asterisque, 96. (1982).

\bibitem{TJ}  Tjurin, A. N; \emph{Classification of n-dimensional vector bundles over an algebraic curve of arbitrary genus},~~~ Izv. Akad. Nauk SSSR Ser. Mat. 30, 1966, 1353--1366.

\end{thebibliography}
\end{document}